\documentclass[letterpaper]{article}

\usepackage{natbib,alifeconf,url,amsmath, tikz, caption, hyperref}  

\usepackage{comment}

\title{Reformalizing the notion of autonomy as closure through category theory as an arrow-first mathematics}

\author{Ryuzo Hirota$^{1}$, Hayato Saigo$^{2}$, \and Shigeru Taguchi$^{3}$   \\
\mbox{}\\
$^1$ Graduate School of Arts and Sciences, University of Tokyo\\
$^2$ Nagahama Institute of Bio-Science and Technology \\
$^3$ Center for Human Nature, Artificial Intelligence, and Neuroscience (CHAIN), \\ Faculty of Humanities and Human Sciences, Hokkaido University \\
hirota@sacral.c.u-tokyo.ac.jp} 

\begin{document}
\maketitle

\begin{abstract}

Life continuously changes its own components and states at each moment through interaction with the external world, while maintaining its own individuality in a cyclical manner. Such a property, known as “autonomy," has been formulated using the mathematical concept of “closure." We introduce a branch of mathematics called “category theory" as an “arrow-first" mathematics, which sees everything as an “arrow," and use it to provide a more comprehensive and concise formalization of the notion of autonomy. More specifically, the concept of “monoid," a category that has only one object, is used to formalize in a simpler and more fundamental way the structure that has been formalized as “operational closure." By doing so, we show that category theory is a framework or “tool of thinking" that frees us from the habits of thinking to which we are prone and allows us to discuss things formally from a more dynamic perspective, and that it should also contribute to our understanding of living systems.
 
\end{abstract}

\section{Introduction}

Life continuously changes its own components and states at each moment through interaction with the external world, while maintaining its own individuality in a cyclical manner. Constructing an artificial system with such a property, i.e., “autonomy," has been one of the major challenges in artificial life \citep[p. 5]{Aguilar2014-kh}.

Mathematical formalization plays a major role in properly capturing such a property. A branch of mathematics called category theory \citep{maclane1971categories, awodey2010category, simmons2011introduction}, which we will introduce later, has been used on occasion in theoretical biology as a formal framework that allows describing relations without relying on concrete components \citep[e.g.,] []{Rosen1991-wg, Varela1979-de, Nomura2007-mi}. And in recent years, there have again been active attempts to construct a system theory based on category theory \citep[e.g.,][]{Capucci2022-gj, Virgo2021-gv, Fong2019-yg}.

In particular, Francisco Varela, who with Humberto Maturana proposed the concept of autopoiesis \citep{Maturana1980-gf, Maturana1987-sh}, also attempted to formalize his theory using concepts of category theory \citep{Varela1979-de, Kauffman2017-wu}. These attempts, however, were eventually regarded as static and were gradually replaced by modeling with dynamical systems \citep[e.g.,][]{Kelso1995-kh, Thelen1994-tp, Di_Paolo2017-dc}, which are supposed to capture dynamic changes in the system. As a result, they have not gained as much attention from the current generation of researchers.

However, category theory can be interpreted as a step forward in the direction of grasping a very dynamic way of being. Although category theory is widely used by modern mathematicians, it is sometimes perceived as a static theory in that it describes the “invariant" structure behind various branches of mathematics. However, if we examine the original intentions of the theory's originators, it becomes clear that category theory was initially developed as a means of capturing the dynamic “movement" inherent in mathematical thinking\footnote{Some readers may think that category theory is static because it does not appear to explicitly incorporate temporal changes, as in dynamical systems. However, as will be mentioned later, (discrete) dynamical systems can be viewed as functors from a monoid to the category of sets \citep[p. 329]{Spivak2014-yx}, and from that perspective, the iteration of time in dynamical systems is nothing other than the composition of the arrows of the monoid. In other words, if dynamical systems are said to be dynamic, then so is a category (especially a monoid). There are also proposed models that allow the structure of a category to change dynamically and stochastically \citep[e.g.,][]{Fuyama2020-rb}.}. In fact, the originators positioned their theory as an extension of the Klein's “Erlangen Program" \citep{klein1893vergleichende} that sought to view geometry as a field centered on dynamic transformations, rather than static shapes to be transformed: “[Category theory] may be regarded as a continuation of the Klein Erlangen Program, in the sense that a geometrical space with its group of transformations is generalized to a category with its algebra of mappings." \citep[][p. 237]{Eilenberg1945GeneralTO}

Such an intention of the originators is deeply and subtly incorporated into the designs of the basic concepts of category theory. In particular, if we look carefully at the “axioms of category," the very basis of category theory, it becomes clear that they are rooted in a perspective that views everything not as a set of point-like elements, but as an “arrow."

Based on this view on category theory, we use category theory not as a general-purpose language that can be applied to the mathematical modeling of general phenomena, but as a “tool of thinking" that leads us to a way of thinking that captures the dynamic nature of things by enabling us to avoid the “habits of thinking" to which we are prone.

As the first small step in this larger endeavor, we attempt here to reformulate a structure called “closure" using category theoretic thinking. It refers to the maintenance of certain characteristics of something constantly changing, and has been adopted by many researchers as a key concept in capturing the autonomous nature of life. One of the examples is the notion of “operational closure," originally proposed by \cite{Maturana1980-gf} and subsequently refined by \cite{Di_Paolo2014-vs}. It describes how the processes that constitute a system are closed with respect to “enabling relation" and, together with the “precariousness" of a network of such processes, is claimed to define the autonomy of living systems and others.

While we believe that this concept is surely effective in capturing the fundamental characteristics of living systems, we also believe that there are some limitations with it.

Given this, from the perspective of category theory as an “arrow-first mathematics," we attempt to shed new light on the autonomous nature of life, which is described as “closure," by means of a structure known as a “monoid," which will be introduced in detail in a later section.

Although the concepts we use in this paper are mathematically very elementary, they nevertheless allow us to formulate in a concise and productive way the discussion on the autonomy of life that has so far been presented in an advanced natural language (or in other mathematical frameworks such as dynamical systems). This shows that category theory is a promising approach for studying complex and dynamic systems, including living organisms.

The structure of the paper is as follows. The next section provides an introduction to the basic concepts of category theory and shows that it incorporates an "arrow-first" perspective. The third section discusses the formalization of the autonomy of life as "closure," focusing on the concept of “operational closure" formally defined by \citep{Di_Paolo2014-vs}. The fourth section introduces the concept of “monoid" (a category with a single object) as the basic mathematical concept to formalize the notion of operational closure in living systems, with a particular focus on their “self-mediating" nature. The discussion will examine the benefits and future prospects of formalization by category theory more generally. Finally, the conclusion will provide a brief summary of the paper's main findings.

\section{Category theory as an “arrow-first" mathematics}

In this section, we introduce the concept of “category" and show that it has an “arrow-first" perspective.

\textit{Definition 1 (Category)}: A category is composed of two kinds of entities, namely, “objects" and “arrows (or morphisms)", that satisfy the following axioms. Any entities and relations that satisfy the axioms can be considered as “objects" and “arrows," respectively, regardless of their specific components.

\textit{Axiom 1 (Arrows and objects)}: Each arrow $f$ has its “domain" (source) object $dom(f)$ and  “codomain" (target) object $cod(f)$. An arrow such that $dom(f)=X$ and $cod(f) =Y$ can be expressed as $f: X \rightarrow Y$ or as follows:

\begin{eqnarray}
\begin{tikzpicture}[auto]
\node (x) at (0, 0) {$X$}; 
\node (y) at (1, 0) {$Y$};
\draw[->] (x) to node {$\scriptstyle f$} (y);
\end{tikzpicture}
\end{eqnarray}

\textit{Axiom 2 (Composition)}: A pair of arrows $f, g$ can be “composed" into $g \circ f$ if the domain of one arrow is equal to the codomain of another, i.e., $cod(f) = dom(g)$. Assuming $f:X \rightarrow Y$ and $g: Y \rightarrow Z$, the following can be expressed:

\begin{eqnarray}
\begin{tikzpicture}[auto]
\node (x) at (0, 0) {$X$}; 
\node (y) at (1.5, 1) {$Y$};
\node (z) at (3, 0) {$Z$};
\draw[->] (x) to node {$\scriptstyle f$} (y);
\draw[->] (y) to node {$\scriptstyle g$} (z);
\draw[->] (x) to node {$\scriptstyle g \circ f$} (z);
\end{tikzpicture}
\end{eqnarray}

\textit{Axiom 3 (Associative law)}: The composition of arrows satisfies the “associative law," i.e.,

\begin{eqnarray}
(h \circ g) \circ f = h \circ (g \circ f)
\end{eqnarray}

This means that assuming $f:X \rightarrow Y$, $g: Y \rightarrow Z$ and $h: Z \rightarrow W$, the following diagram is “commutative," i.e., no matter which path the arrows are composed through, if the start and end points are the same, the result is the same:

\begin{eqnarray}
\begin{tikzpicture}[auto]
\node (x) at (0, 0) {$X$}; 
\node (y) at (1.5, 0) {$Y$};
\node (z) at (1.5, 1.5) {$Z$};
\node (w) at (3, 1.5) {$W$};
\draw[->] (x) to node {$\scriptstyle f$} (y);
\draw[->] (y) to node {$\scriptstyle g$} (z);
\draw[->] (z) to node {$\scriptstyle h$} (w);
\draw[->] (x) to node {$\scriptstyle g \circ f$} (z);
\draw[->] (y) to node[swap] {$\scriptstyle h \circ g$} (w);
\end{tikzpicture}
\end{eqnarray}

\textit{Axiom 4 (Identity)}: Each object $X$ has its corrresponding arrow to itself $1_X: X \rightarrow X$ called “identity" such that for any arrow $f: X \rightarrow Y$,

\begin{eqnarray}
f \circ 1_X = 1_Y \circ f = f
\end{eqnarray}

Intuitively, we can think of an object as representing a “thing," “phenomenon," or “event," and an arrow as representing a directed “relation," “process," or “transformation" between them. The most typical example that appears in the natural sciences is a category in which the distinct states of a system are considered as the objects and possible transitions between the states as the arrows (\cite{Saigo2019-on} name this category the “category of mobility").

According to the basic idea of category theory, it does not begin by assuming objects unrelated to arrows; rather, an object is characterized only by what kind of arrows it has to other objects (and to itself). In this context, it is important to note that arrows are not reduced to pairs of objects. In general, there can be more than one arrow with the same domain and codomain, as schematically depicted in the diagram below: 

\begin{eqnarray}
\label{diagram1}
\begin{tikzpicture}[auto]
\node (A) at (1, 1) {$X$}; 
\node (B) at (3, 1) {$Y$};
\node (X) at (0, 1) {};
\node (Y) at (0, 1.5) {};
\node (Z) at (0, 0.5) {};
\node (F) at (4, 1) {};
\node (G) at (4, 1.5) {};
\node (H) at (4, 0.5) {};
\draw[<-] (X) to (A); \draw[->] (Y) to (A); \draw[->] (Z) to (A);
\draw[<-] (F) to (B); \draw[->] (G) to (B); \draw[->] (H) to (B);
\draw[->] (A) to (B);
\draw[->] (A) to[bend left=30] (B);
\draw[->] (A) to[bend left=60] (B);
\draw[->] (A) to[bend right=30] (B);
\draw[->] (A) to[bend right=60] (B);
\end{tikzpicture}
\end{eqnarray}
Therefore, if we obtain only information about the objects from a category, we cannot recover the original category, whereas if we are given information about the arrows, we can reconstruct a complete picture of the original category. Unlike the set-theoretic approach, which typically first considers objects (sets) and then discusses the relations (functions) between them, it is precisely arrows that play the leading role in category theory.

Also important is the “identity law," which allows us to identify an object with its corresponding “identity" arrow and thus to regard an object as a kind of arrow. In other words, since a category consists of objects and arrows, and the objects are also a kind of arrow, we can say that a category is actually composed of various kinds of arrows. In this sense, again, “It's the arrows that really matter!" \citep[][p. 8]{awodey2010category}

Based on these considerations, we can interpret category theory as an \textit{“arrow-first"} mathematics, where everything is conceived of as an arrow. In this perspective, objects, which may appear static and fixed, are considered a particular type of arrow that are dedicated to \textit{mediating} between arrows. Thus, category theory as “arrow-first mathematics" considers everything as an “arrow," while acknowledging the need for the identity arrows, which correspond to the objects. As a result, category theory does not merely reduce everything to transitive relations but rather always includes something individual and intransitive. This makes category theory different from simplistic relationalism or relational monism, which dissolve individuals into relations.\footnote{This stance of category theory is closely aligned with the one expressed by \cite{Varela1976-sa,Varela1979-de} in his phrase “Not one, not two." According to it, the dyad of “the it" and “the process leading to it" can be understood neither by separating them dualistically nor by reducing one to the other monistically; they can only be understood as complementary. Moreover, it is emphasized that they are not in a symmetrical relationship where the two are mutually exclusive, but in an asymmetrical relationship where one emerges from the other. Although it is beyond the scope of this paper, Varela analyzed such a complementary relationship using the concept of “adjunction" in category theory \citep[pp. 97-99]{Varela1979-de}.}

Furthermore, as mentioned above, there can be many (sometimes infinite) arrows between two objects in general (as depicted in the diagram (\ref{diagram1})). Typically, we tend to think of a “relation" as being uniquely determined by the pair of objects related. However, in category theory, arrows are not reduced to pairs of objects. Therefore, the relations represented by arrows are not reduced to elements but are inherently diverse and pluralistic. These relations are more flexible than what we typically associate with the term “relation," and it may be more appropriate to use the term “mediation," which we will discuss in detail later.

We can also consider a category whose objects are arrows in another category, as we will discuss later. In that sense, what constitutes an arrow and an object is by no means predetermined. Rather, it depends on the context in which they are being used. 

In summary, category theory characterizes mathematical objects not by what they are composed of, but by how they behave and relate to other things, i.e., what kind of arrows they have to and from other objects. Moreover, it views even the objects themselves as a form of arrow or process of “being an object," considering it an indispensable aspect. This nature of category theory would be compatible with artificial life, which seeks to characterize life by its relational and behavioral properties rather than by specific components such as proteins or DNA, but also to implement it as an individual.

In the following sections, we will consider the characteristics of living systems through the lens of category theory as an “arrow-first mathematics," which allows us to grasp the dynamic nature of things, as described in this section. In particular, our focus will be on the concept of “closure," which we will formalize using a special case of category called “monoid."

\section{Autonomy and Closure}

In theoretical biology, the notion of “closure" has often been used to capture the distinctive characteristics of life, particularly its autonomous nature \citep{Moreno2015-cd}. Simply put, the notion of closure refers to the phenomenon of returning to the original state after some manipulations, processes, or changes in a specific sense, i.e., being “closed" with respect to them. The notion of closure has been considered in various forms in theoretical biology. For instance, \cite{Rosen1991-wg} pointed out that the metabolic processes of life are “closed" with respect to “efficient causation" and attempted to formalize them using category theory \citep[see also][]{Letelier2003-cm}. Furthermore, \cite{Montevil2015-su} argue that the constraints on processes in living systems are closed: a process constrained by the outcome of another process creates a constraint on yet another process, forming a cycle of constraints.

Here, we particularly focus on what is known as “operational closure." It basically means that some processes are “closed" with respect to their operations, and was originally proposed by \cite{Maturana1980-gf} as a property of “autopoietic" systems. It was also developed by \cite{Varela1979-de} to characterize autonomous systems in general, not limited to biochemical processes in cells (including nervous, immune, or social systems). Subsequently, \cite{Di_Paolo2014-vs} define it more formally as the property that every process constituting a system is enabled by at least one other constituting process and also enables at least one other constituting process. In other words, it refers to how the processes that constitute the system are “closed" with respect to “enabling relations." 

The concept of operational closure has become a foundation for the “enactive approach" in cognitive science, which argues that autonomy is at the basis of cognition. In contrast to mainstream cognitive science's computational and representationalist premises, the enactive approach, proposed by \cite{Varela1991-xr}, highlights that the agent and the world are not independently given but are brought forth or “enacted" through their interactions \citep{Thompson2007-po, Di_Paolo2017-dc, Di_Paolo2018-qd}. The enactive approach is built on the theoretical traditions of autopoiesis and autonomy, but it also attempts to modify and expand them in fundamental ways \citep{Weber2002-hf, Di_Paolo2005-mw, Di_Paolo2009-ze, Froese2010-ub}. The redefinition of operational closure proposed by \cite{Di_Paolo2014-vs} is part of this effort to refine and extend the theoretical framework of autopoiesis and autonomy.

\begin{figure}[t]
\begin{center}
\includegraphics[width=2.4in]{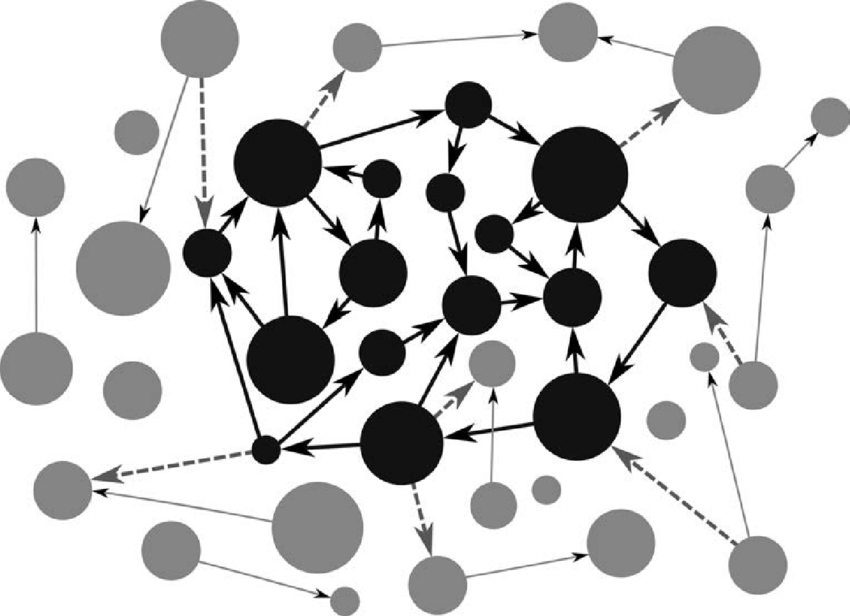}
\caption{A schematic illustration of the concept of operational closure \citep[reprinted from][]{Di_Paolo2014-vs}. Copyright Ezequiel Di Paolo, 2013. This work is licensed under a \href{http://creativecommons.org/licenses/by-nc-sa/3.0/deed.en_US}{Creative Commons Attribution-NonCommercial-ShareAlike 3.0 Unported License}.}
\label{fig3}
\end{center}
\end{figure}

Figure~\ref{fig3} represents the concept of operational closure \citep{Di_Paolo2014-vs} . The diagram shows a series of circles that represent processes, and arrows that represent enabling relations between processes. When one process cannot occur without another, there is an enabling relation from the latter to the former \citep{De_Jaegher2010-md}. The circles in black represent processes that are closed with respect to enabling relations, meaning that each black circle is enabling and enabled by another process represented by a black circle. The processes represented by these black circles, taken together, constitute an operationally closed system.

\begin{figure}[t]
\begin{center}
\includegraphics[width=2.4in]{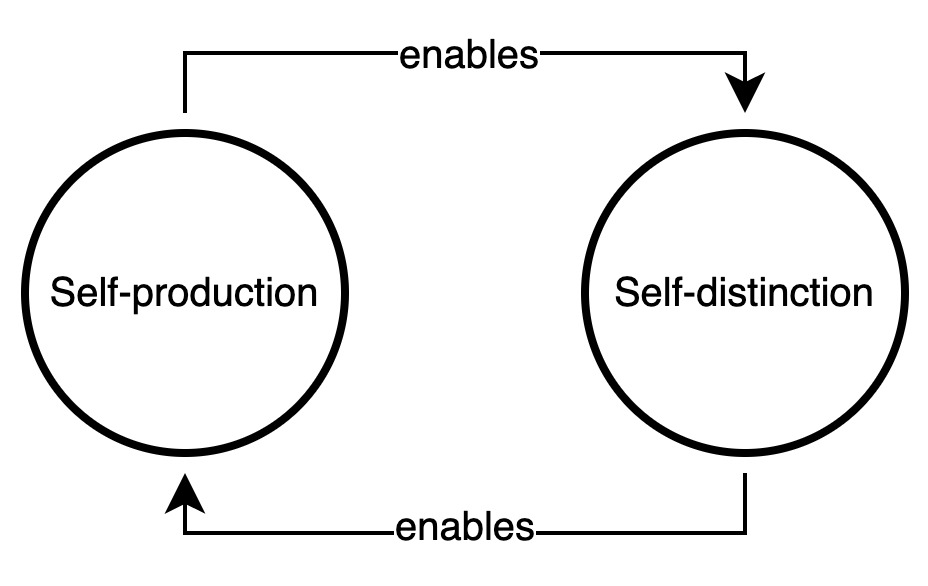}
\caption{Operational closure between self-production and self-distinction.}
\label{fig4}
\end{center}
\end{figure}

Typical examples of operationally closed processes are “self-distinction" and “self-production" in a living cell (Figure~\ref{fig4}). The formation of a membrane creates a distinct physicochemical condition within it that facilitates metabolic reactions (self-distinction), while such metabolic reactions produce the components necessary for membrane formation (self-production), leading to a mutually enabling relation known as “autopoiesis," which was proposed by Maturana and Varela \citep{Maturana1980-gf, Maturana1987-sh, Varela1997-vn, Di_Paolo2018-fw}.

Importantly, \cite{Di_Paolo2014-vs} argue that for a system to be autonomous, it is not enough that it is operationally closed; it must also be additionally “precarious." Precariousness means that the operationally closed processes are not given, but are non-trivial, produced and maintained by themselves. It is formally defined as the property that if any of the enabling relations cease to hold, the entire system cannot exist \citep{Di_Paolo2014-vs}. 

While the concept of operational closure is surely effective in capturing the autonomous nature of living systems, we consider that there are also some limitations with it. First of all, the enabling relation between processes (which are themselves a kind of relation) is so abstract and intangible that we cannot easily imagine, and as a result of trying to visualize such a hard-to-imagine thing, the visualization can easily cause a kind of reification of those processes. Of course, the authors emphasize that what are represented by the circles in the figure are not thing-like entities but processes, however, in our natural way of thinking, it is hard to imagine them properly, let alone their relationships. Therefore, there is a need for a concrete framework that allows for a more tangible understanding of these abstract relationships, other than arguments in natural language and/or illustrations, as \cite{Beer2020-by} also points out. As noted by \cite{Di_Paolo2022-fl}, concepts in theories of autopoiesis and enaction are often misinterpreted and misapplied by researchers from outside the field, and we believe that this is partly due to the abstract nature of these concepts.

Furthermore, the figure could induce a way of thinking in which the processes depicted by the circles are viewed as spatially and/or functionally differentiated and modularized, similar to the schematic diagram of computational architecture. While this is not the intended interpretation of the authors, it is evident that the figure could be associated with functionalist and computationalist perspectives, which enactivists are supposed to be avoiding \citep{Di_Paolo2017-dc, Di_Paolo2018-qd}.

To address these challenges, in the next section, we will attempt to reformulate the concept from a more fundamental and simple standpoint by using the language of category theory introduced in the previous section. Category theory is precisely a branch of mathematics that explores the properties of relationships between objects and makes it possible to concretize even highly advanced concepts, such as “relations of relations (of relations)," which are difficult (or almost impossible) to capture by natural language discussion alone. Our aim is not to refute the existing definition and depiction but to expand their potential by shedding light on them from a different perspective.

Also, the authors argue that “although the choice of processes under study is more or less arbitrary and subject to the observer’s history, goals, tools, and methods, the topological property unraveled is not arbitrary" \citep[][p. 71]{Di_Paolo2014-vs}, meaning that the definition's essence lies not in the specific processes represented by the circles (which vary depending on the observer) but rather in the properties of the enabling relations among them. As described in the last section, category theory is the mathematical framework created precisely to allow us to discuss such properties of relationships among objects without relying on their specific components. In this sense, too, our approach aligns with and complements the authors' intent.

\section{Mediation and Monoid structure}

For the purpose described in the last section, we will use the simplest mathematical formalization of the structure of what is called a “closure" in general: namely, a “monoid," a category with only one object. Importantly, even if it has only one object, there can be an infinite number of arrows (this follows from the irreducibility of an arrow to the pair of objects). This property allows us to successfully capture the way something is always varied, yet remains the same in a specific sense.

\textit{Definition 2 (Monoid)}: A monoid is a category with a single object. 

Since it has only one object, a monoid is, in effect, a “collection of arrows." And since the domain and codomain of all the arrows are the same, it is possible to freely composite the arrows with each other. A monoid can be expressed as follows:

\begin{eqnarray}
\begin{tikzpicture}[auto]
\node (A) at (1, 1) {$\cdot$}; 
\draw [->] (A) edge[loop above] node{} (A);
\draw [->] (A) edge[loop below] node{} (A);
\draw [->] (A) edge[loop right] node{} (A);
\draw [->] (A) edge[loop left] node{} (A);

\end{tikzpicture}
\end{eqnarray}

The concept of operational closure can be formulated as a monoid, insofar as it exhibits the characteristics of a “closure." In the following, we aim to achieve this by specifically focusing on “enabling relations" as arrows, following the basic stance of category theory as “arrow-first mathematics."

What is significant about an enabling relation is that it is not deterministic; “A enables B" does not mean that B will certainly occur or exist if A is present. Rather, what it exactly means is that \textit{without A, there would be no B} \citep{De_Jaegher2010-md}. This difference, often described as the one between “determination" and “dependence" \citep[e.g.,][p. 337]{Di_Paolo2018-qd}, is crucial. Consider, for example, the relationship between seed and germination. Since germination is affected by various factors including environmental ones, it cannot be said that just because a seed exists, it will necessarily sprout. However, this does not mean there is no law and everything is uncertain. Rather, it is quite certain that there is a particular relationship: “If there is no seed, there will be no sprouting."  This kind of relationship can be found everywhere in life phenomena in general. For example, one of the characteristics of life, at least on Earth today, is that individuals do not spontaneously arise, or simply put, “no children without parents" or “life only comes from life" \citep{Oono2012-jn, Froese2019-nj}. In this context, the relationship of “without A there would be no B" captures the nature of life as being inherently path-dependent and historical \citep[see also][]{Longo2012-lw}.

More generally, what has been studied as “causality" in biology, neuroscience, and medicine is essentially this “without A, there is no B" relation\footnote{Interestingly, such a conception of causality as a fundamental dependency on others can be closely related to what is called “dependent arising" (paṭicca samuppāda) in the Buddhist tradition, to which \cite{Varela1991-xr} were also profoundly concerned.}: the relation between genetic “knockouts" and their effects on phenotype, the relation between physical and physiological changes in the nervous system and the transformation of subjective experience, and so on. Furthermore, a similar view of causality based on counterfactuals can also be found in the framework of statistical causal inference proposed by Pearl and colleagues \citep{Pearl2018-yz}, allowing for quantitative as well as qualitative analysis.

The relationship of “without A, there would be no B" can be conceptualized as one of “mediation" \citep{Taguchi2019-fj}. This is exemplified by neural processes in the brain mediating human behavior or honeybees mediating the pollination of flowers, both of which imply that the former is necessary for the latter to occur. Moreover, the relationship “without A, there would be no B" does not exclude that it involves other factors than A and B. In essence, the relationship “without A, there would be no B" differs from the deterministic relationship “if there is A, there is always B" in that it typically involves other variables besides A that contribute to the occurrence or existence of B. For instance, the mere presence of a seed does not guarantee it will sprout; additional environmental factors such as soil, water, and heat are necessary.

Based on the notion of mediation as such, we can further obtain the perspective that all things are mediated by each other and things can exist only through such various forms of mediation. In other words, rather than something unmediated existing first and then entering into a mediating relationship, \textit{mediation always comes first}, and what seems unmediated emerges only as a \textit{mediating} point between mediating relations.

This mediation-based perspective is in deep accordance with the “arrow-first" perspective of category theory described above, in which objects do not exist independently as themselves but are characterized only as the “hubs" of arrows. In other words, the concept of category can accurately represent a worldview that sees everything as mediation.

From the perspective discussed so far, let us again consider the nature of living systems, especially the autonomy (and selfhood) expressed as operatinally closed, using category theory as a “tool of thinking".

In our lives, there are innumerable mediating relations among things, and some of them return to the original in a cyclic manner. They include horizontal relations with the outside, such as “agent → environment → agent" or “agent → other agents → agent," as well as vertical relations between the global and local inside the agent, such as “organism → organ → organism." They have been referred to as “circular (or reciprocal) causality" in the enactive approach and elsewhere \citep[e.g.,][]{Thompson2007-po, Fuchs2017-pm, Fuchs2020-pc, Tschacher2007-lk}.

One typical example of this is chemotaxis. It refers to the tendency of microorganisms like E. coli to self-migrate towards environments richer in nutrients. In this case, the agent's existence as mobile allows for the presence of a specific environment around it through self-movement, and in turn, the environment enables the agent to survive. Thus, there is a mutually enabling and mediating relationship between the agent and the environment around it. In other cases, known as “niche construction" \citep{Odling-Smee2003-uf} or “extended physiology" \citep{Turner2000-mo}, the agent “mediates" the environment more directly. As an example, Di Paolo (2009) describes how the water boatman, an aquatic insect, is able to keep air bubbles on its body surface underwater by means of hairs with water-repelling properties, which allows it to breathe and spend more time underwater, creating a mutually mediating relationship between the water boatman and its environment: water boatman → air bubbles → water boatman.

Considering a category whose arrows are the mediating relations between various objects, and then focusing on a single object and its various arrows from itself to itself (i.e., self-mediation), a monoid can be obtained (as a subcategory). As noted above, in category theory, relations represented by arrows are not reducible to the pairs of objects since the arrows between objects are more than one in general, as illustrated in diagram(\ref{diagram1}). In particular, there can be multiple or even innumerable arrows from an object to itself, with the “identity" arrow being the most trivial example. In the monoid considered here, the most trivial arrow, the “identity," is the tautological self-mediation “self → self," i.e., “without the self, there would be no self." However, as mentioned earlier, in category theory, an object itself does not have any characteristics and is characterized only in terms of its relations to others, i.e., arrows. This is also true of the monoid we are considering here, meaning that the self as the object does not exist independently of its relations to others, but can only exist as a “hub" through which the arrows are connected to each other. In other words, the “self" as the object of the monoid, or the self-mediation “self → self" as the identity arrow, can only exist through other forms of self-mediation, such as “self → environment → self," “self → other → self," or “self → organ → self."

\cite{Varela1991-dc} states that the autonomous self is “a meshwork of selfless selves" that is interwoven with various processes, including metabolic, immune, sensorimotor, and social ones, with no single substantial core. This conception of the self can be expressed in a natural way by the formalization of the self as a monoid described above.

Such a depiction of the autonomous self as a monoid might be interpreted as solipsistic, but this is not the case. In fact, the arrows of a monoid can be “factorized," which allows us to recover other factors as regularities that appear among the arrows. In the first place, an isolated factor cannot constitute mediation; mediation means that the existence of a factor intrinsically involves the existence of other factors. Therefore, the picture in which the self as an object can be grasped only in terms of its various mediations is thoroughly “world-involving" \citep{Di_Paolo2017-dc}. 

Furthermore, category theory's “arrow-first" perspective stands in contrast to simplistic relational monism, which reduces everything into a web of relations and denies individuality. This aligns with the stance of the enactive approach, particularly of the one referred to as “autopoietic enactivism," which views the self as emerging from/through environmental interactions rather than as being a pregiven entity, yet without entirely dismissing the notion of an autonomous self as an illusion or observer's artifact \citep{Barandiaran2017-rp, Di_Paolo2009-ze}, as expressed in the phrase “neither individualistic, nor interactionist" in the context of sociality \citep{Di_Paolo2017-ku}. Thus, our view of autonomy as the monoid can also address the concern raised by \citet{Barandiaran2017-rp} that autonomy and the coupling with the environment are likely to be seen as a binary choice or tradeoff; What essentially constitutes the monoid are the arrows that represent self-enabling involving the environment, but these arrows are mediated by one object, namely the autonomous self, which in turn would be meaningless if isolated from the arrows.

One merit of using category theory as a language for formalization is that it allows for constructing a category from another category. One such construction is a “coslice category,"  whose object is an arrow in another category. More concretely, the objects and arrows of a coslice category $A/\Omega$ of a category $\Omega$ (where $A$ is an object arbitrarily selected from all the objects of category $\Omega$) are defined as follows:

\textit{Objects in $A/\Omega$}: Arrows in category $\Omega$ with $A$ as its domain, such as $f: A \rightarrow X$ and $g: A \rightarrow Y$.

\textit{Arrows in $A/\Omega$}: An arrow in $A/\Omega$ between objects $f: A \rightarrow X$ and $g: A \rightarrow Y$ is defined as a triplet $(f, g, t)$ such that an arrow $t: X \rightarrow Y$ in category $\Omega$ satisfies $g = t \circ f$, i.e., the following diagram is commutative in category $\Omega$:

\begin{eqnarray}
\begin{tikzpicture}[auto]
\node (A) at (0.75, 0) {$A$}; 
\node (X) at (0, 1) {$X$}; 
\node (Y) at (1.5, 1) {$Y$}; 

\draw[->] (A) to node {$\scriptstyle f$} (X);
\draw[->] (A) to node[swap] {$\scriptstyle g$} (Y);
\draw[->] (X) to node {$\scriptstyle t$} (Y);

\end{tikzpicture}
\end{eqnarray}

\begin{figure}[htbp]
\begin{center}
    \begin{tabular}{cc}
        \begin{minipage}[t]{0.4\linewidth}
            \begin{center}
            \begin{tikzpicture}[auto]
                \node (A) at (1.2, 1) {$A$}; 
                \draw [->] (A) edge[loop above] node{$\scriptstyle f_1$} (A);
                \draw [->] (A) edge[loop below] node{$\scriptstyle f_3$} (A);
                \draw [->] (A) edge[loop right] node{$\scriptstyle f_2$} (A);
                \draw [->] (A) edge[loop left] node{$\scriptstyle 1_A$} (A);   
            \end{tikzpicture}
            
            \end{center}
        \end{minipage} &
        \begin{minipage}[t]{0.4\linewidth}
            \begin{center}
            \begin{tikzpicture}
                \node (1_A) at (1, 1) {$1_A$}; 
                \node (f1) at (2, 2) {$f_1$}; 
                \node (f2) at (3, 1) {$f_2$}; 
                \node (f3) at (2, 0) {$f_3$}; 
                \draw[->] (1_A) to node {} (f1);
                \draw[->] (1_A) to node {} (f2);
                \draw[->] (1_A) to node {} (f3);
                \draw[->] (f1) to node {} (f2);
                \draw[->] (f2) to node {} (f3);
            \end{tikzpicture}
            
            \end{center}
        \end{minipage}
    \end{tabular}
\end{center}

\caption{Monoid $\Omega$ (left) and Coslice category $A/\Omega$ (right).}
\label{fig5}
\end{figure}
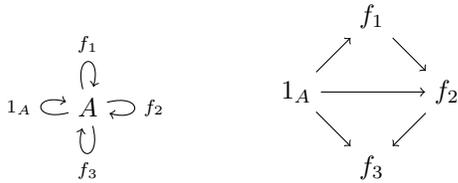

Considering the coslice category of the monoid we have considered in this section, it is possible to take as objects the mediating relations that are the arrows in the monoid and to explicitly explore the relations between them (Figure~\ref{fig5}). This is probably close to the figure of operational closure in \cite{Di_Paolo2014-vs} (Figure~\ref{fig3}), in which each process is differentiated and depicted as an individualized circle.

For example, the mutual enabling relationship between self-distinction and self-production (illustrated in Figure~\ref{fig4}) may be associated with the structure in the coslice category of the monoid. 

The identity arrow $1_A$ of the object $A$ in the original category $\Omega$ is included as well in the coslice category $A/\Omega$ as one of the objects and is placed along with other objects in the coslice category (i.e., arrows in the original category). However, this object is a special one known as the “initial object," which has an arrow to every object in the category, since any arrow $f: A \rightarrow X$ in the original category $\Omega$ and the identity arrow $1_A$ make the diagram below commutative and thus there is always an arrow in $A/\Omega$ from $1_A$ to $f$: 

\begin{eqnarray}
\begin{tikzpicture}[auto]
\node (A) at (0.75, 0) {$A$}; 
\node (X) at (1.5, 1) {$X$}; 
\node (A2) at (0, 1) {$A$}; 

\draw[->] (A) to node[swap] {$\scriptstyle f \circ 1_A = f$} (X);
\draw[->] (A) to node {$\scriptstyle 1_A$} (A2);
\draw[->] (A2) to node {$\scriptstyle f $} (X);
\end{tikzpicture}
\end{eqnarray}
Such a property of an object possessing a unique arrow to each of all other objects in a category is known as the “universal property," and it plays a major role in the definition of various crucial concepts in category theory, such as those of “product" and “equalizer" \citep{Leinster2014-cy}; These concepts are defined as the distinctiveness of a specific object in a category in terms of its relation with all other objects.

Back to the discussion on autonomy, the concept of universal property may also play an important role in characterizing the particularity of self-distinction in an autonomous system. According to \cite{Virgo2011-cy}, self-distinction, i.e., being a distinct unity (with spatial boundaries) is, on the one hand, nothing more than one of the constitutive processes (as depicted in Figure~\ref{fig4}), but, on the other hand, it is still an exceptional one in that “it enables a great number of processes" (p. 247). 
In other words, self-distinction, although seemingly static, can be thought of as one, yet special, kind of process that mediates various, or perhaps even all, other self-producing processes. This would be depicted as a black circle with an arrow to all other black circles in a diagram as in Figure~\ref{fig3}. 

This feature of self-distinction may be associated with the universal property of the identity arrow $1_A$ (representing “being an object" as a process) of monoid $\Omega$ as the initial object of the coslice category $A/\Omega$, possessing arrows to all other objects (Figure~\ref{fig5}).This category-theoretic understanding of self-distinction, focusing on temporal and processual properties rather than spatial ones, can provide new insights into artificial life by generalizing the role of membranes.

Our intention in this section is not to argue that anything described as a monoid is autonomous or living. On the contrary, anything possessing a recursive property can be described as a monoid. This point would be partly related to the argument that autonomy requires not only being operationally closed but also being “precarious," that is, if any of the enabling relations between the constitutive processes disappear, the entire system ceases to exist \citep{Di_Paolo2014-vs}. Such properties that distinguish the living from the non-living might be described as a specific interdependent structure in the coslice category, which needs to be addressed in future studies.

\section{Discussion}

Finally, in this section, we explicitly highlight some of the merits of employing category theory as a language for formalization.

First of all, category theory allows us to deal with the concept of “relation," which is so multifaceted that it often obscures the discussion, using a more precise concept: arrow. The word “relation" can sometimes be interpreted as reducible to the two terms it relates. In contrast, as discussed in this paper, arrows in category theory represent a specific and enriched notion, which we refer to as “mediation."

Another advantage of category theory as a language for formalization is that it allows us to associate between structures across different hierarchies. For example, it is possible to think of a part of an object (e.g., an organ of the self) or even what includes the object as a part (e.g., society to which the self belongs) as another object. In general, category theory can treat the part and the whole (or the local and the global) as being equal; they are equally treated as objects, and the relationship of containment between them is represented as an arrow between objects. Unlike the set-theoretic approach, which is based on and privileges the containment relationship “something is an element of another thing," category theory allows from the beginning to deal with relations in a broader sense as its default. This nature of category theory can make it possible to speak consistently about the hierarchical, vertical (local-global) relationship within the system and the horizontal relationship between the agent and the environment (including other agents). 

In addition, although not discussed in this paper, category theory also allows us to rigorously describe higher-order relations, such as a correspondence between categories (an arrow in a category in which each object is a category), which is called a “functor," and even a consistent correspondence between functors (an arrow in a category in which each object is a functor), which is called a “natural transformation." It would be nearly impossible to speak of such higher-order relationships strictly in natural language alone. Category theory, in contrast, was designed by the originators from the beginning to address these higher-order relationships. As \citet{Leinster2014-cy} notes, “In fact, it was the desire to formalize the notion of natural transformation that led to the birth of category theory" (p. 9). And even in such categories as the “category of categories" and the “category of functors," the logic about autonomy discussed in this paper can hold true. Hence, these concepts could be applied to, for example, the formalization of the autonomy at the higher-order level discussed in the enactive approach, such as the autonomy of the sensory-motor schemes \citep{Di_Paolo2017-dc} and the autonomy of social interactions \citep{Di_Paolo2018-qd}.

Finally, it should also be noted that dynamical systems, which have been frequently used in theoretical biology and cognitive science, including the enactive approach \citep[e.g.,][]{Kelso1995-kh, Thelen1994-tp, Clark1998-xj,  Beer2004-us, Di_Paolo2017-dc} are one of the typical examples of monoids. The category of the “idea" of dynamical systems is a monoid whose arrows are generated by $n$-times compositions ($n = 1, 2, ...$) of a single arrow, and each discrete dynamical system can be regarded as a “set-valued functor" from this category to the category of sets (a category whose objects and arrows are sets and mappings between them). Thus, monoids can be seen as a relaxation of the deterministic property of dynamical systems, in which all arrows are generated from a single arrow. From our perspective, it is the structure of “repetition" that people have sought to capture using the language of dynamical systems, and the commitment to determinism is not necessarily essential. 
Indeed, the enactive approach is beginning to question the assumption of structural determinism \citep{Froese2019-nj, Froese2023-ej, Fuchs2021-ef} that was central to Maturana and Varela's original concept of autopoiesis \citep{Maturana1980-gf, Froese2010-ub}. In this context, monoids, and category theory in general, potentially serve as a more suitable mathematical framework than dynamical systems.

\section{Conclusion}

In this paper, we have introduced category theory as an “arrow-first" mathematics to provide a comprehensive and concise formalization of the autonomous nature of life, especially the structure that has been expressed as “closure." Along the way, category theory has served not merely as a “neutral" formal language but as a tool for thinking that frees our minds from the habits into which they often fall, and leads us to view everything as a relation, or “mediation." This allowed us to formalize the notion of autonomy as a monoid closed under mediating relationships, and to see the autonomous self as a “hub" through which various self-mediating processes are mediated.

Artificial life has become one of the most interdisciplinary and hybrid fields, serving as the nexus of various fields such as biology, chemistry, artificial intelligence, robotics, cognitive science, and philosophy, which makes itself a very “mediating" field of research. Category theory shall serve as a promising platform in such a field for exchanging and integrating ideas from various areas in a comprehensive and productive manner and guiding our thinking in a freer direction.

\section{Acknowledgements}
This work was partially supported by JSPS KAKENHI Grant Number JP20H00001 and JP23H04831. We would like to thank Tom Froese for his insightful discussion on an earlier version of this work and anonymous reviewers for their helpful comments.

\footnotesize
\bibliographystyle{apalike}
\bibliography{reference} 

\end{document}